\documentclass[journal]{IEEEtran}
\usepackage{cite}
\usepackage{times,amsmath,amssymb,psfrag}
\usepackage{graphicx}
\usepackage{url}
\usepackage{bm}
\usepackage{multicol}
\graphicspath{{./pics/}} \DeclareGraphicsExtensions{.pdf,.jpg,.png}
\usepackage{color}
\usepackage{hyperref}
\usepackage{epsfig}
\usepackage{float}
\usepackage{array}
\usepackage{multirow}
\usepackage{etoolbox}
\usepackage{nomencl}
\usepackage{algorithm}
\usepackage{algorithmicx}
\usepackage{algpseudocode}
\usepackage{amsmath}

\newcommand\Nom[3][X]{\nomenclature[#1#3]{#2}{#3}}
\renewcommand\nomgroup[1]{%
  \item[\normalsize\itshape\bfseries
  \ifstrequal{#1}{I}{Indices and Sets}{%
  \ifstrequal{#1}{P}{Notation for MDP-based Model}{%
  \ifstrequal{#1}{N}{Notation for Solution Method}{%
  \ifstrequal{#1}{X}{Other Symbols}{}}}}]%
}

\makenomenclature
\begin{document}
%\begin{multicols}{2}
%\title{MDP-based Resilience Enhancement for Distribution Systems: An Approximate Dynamic Programming Approach}
\title{Markov Decision Process-based Resilience Enhancement for Distribution Systems: An Approximate Dynamic Programming Approach}
\author{
Chong Wang,~\IEEEmembership{Member,~IEEE}, Ping Ju,~\IEEEmembership{Senior Member,~IEEE}, Shunbo Lei,~\IEEEmembership{Member,~IEEE}, \\ Zhaoyu Wang,~\IEEEmembership{Member,~IEEE}, Yunhe Hou,~\IEEEmembership{Senior Member,~IEEE}

\thanks{
	
%This work was supported by the U.S. Department of Energy Office of Electricity Delivery and Energy Reliability.
	
%C. Wang and P. Ju are with the College of Energy and Electrical Engineering, Hohai University, Nanjing 211100, China (e-mail: chongwang@hhu.edu.cn, pju@hhu.edu.cn).

%S. Lei and Y. Hou are with the Department of Electrical and Electronic Engineering, The University of Hong Kong, Hong Kong (e-mail: leishunbo@eee.hku.hk, yhhou@eee.hku.hk).

%Z. Wang is with the Department of Electrical and Computer Engineering,Iowa State University, Ames, IA 50011, USA (email:wzy@iastate.edu).

%Chong Wang is with the Department of Electrical and Computer Engineering,Iowa State University, Ames, IA 50011, USA (email:wangc@iastate.edu).

%Jianhui Wang and Dongbo Zhao are with Argonne National Laboratory, Argonne, IL 60439 USA (email: jianhui.wang@anl.gov, dongbo.zhao@anl.gov). 
}
}

\date{}
\maketitle
%\end{multicols}
%\vspace{-10 mm}
\begin{abstract}
Because failures in distribution systems caused by extreme weather events directly result in consumers' outages, this paper proposes a state-based decision-making model with the objective of mitigating loss of load to improve the distribution system resilience throughout the unfolding events. The sequentially uncertain system states, e.g., feeder line on/off states, driven by the unfolding events are modeled as Markov states, and the probabilities from one Markov state to another Markov state throughout the unfolding events are determined by the component failure caused by the unfolding events. A recursive optimization model based on Markov decision processes (MDP) is developed to make state-based actions, i.e., system reconfiguration, at each decision time. To overcome the curse of dimensionality caused by enormous states and actions, an approximate dynamic programming (ADP) approach based on post-decision states and iteration  is used to solve the proposed MDP-based model. IEEE 33-bus system and IEEE 123-bus system are used to validate the proposed model. 
\end{abstract}
\IEEEpeerreviewmaketitle
\printnomenclature[2.2cm]

\Nom[I01]{$c$}{Index of components.}
\Nom[I02]{$l$}{Index of lines.}
\Nom[I03]{$k,k'$}{Index of terminal buses of line $l$.}
\Nom[I04]{$t, \tau$}{Index of time periods.}
\Nom[I05]{$i,i',j,j'$}{Index of states.}

\Nom[I06]{$\mathcal{A}$}{Set of actions.}
\Nom[I07]{$\mathcal{B}$}{Set of buses.}
\Nom[I08]{${\mathcal{\bar B}}_{i,t}$}{Set of non-islanded buses under the state $S_{i,t}$.}
\Nom[I09]{$\mathcal{\tilde B}$}{Set of substation nodes.}
\Nom[I10]{${\mathcal{C}^{f}_{i,t}}$}{Set of repaired components under state $S_{i,t}$.}
\Nom[I11]{$\mathcal{F}_t$}{Set of all possible failure components at $t$.}
\Nom[I12]{$\tilde{\mathcal{F}}_t$}{Set of actual failure components at $t$.}
\Nom[I13]{${\mathcal{L}}_{i,t}$}{Set of non-islanded lines under state $S_{i,t}$.}
\Nom[I14]{${\mathcal{L}}_{i',t}$}{Set of dispatchable lines under post-decision state $S^{a_t}_{i',t}$.}
\Nom[I15]{${\mathcal{L}}^{nd}_{i,t}$}{Set of non-dispatchable lines under state $S_{i,t}$.}
\Nom[I15]{${\mathcal{L}}^{d}_{i,t}$}{Set of dispatchable lines under state $S_{i,t}$.}
\Nom[I16]{${\mathcal{N}}_{k,i,t}$}{Set of nodes connected to bus $k$ under state $S_{i,t}$.}
\Nom[I17]{$\tilde{\mathcal{R}}_{\tau}$}{Set of repaired components at $\tau$.}
\Nom[I18]{${\mathcal{S}}$}{Set of Markov states.}
\Nom[I18]{${\mathcal{S}_i^{post}}$}{Set of post-decision states of state $S_{i,t}$.}
\Nom[I19]{${\mathcal{T}}$}{Set of time periods.}

\Nom[P01]{$a_t$}{Action at $t$.}%
\Nom[P02]{$C_t$}{Immediate cost at $t$.}%
\Nom[P02]{$C_l$}{Operational cost of line $l$ at $t$.}%
\Nom[P03]{$F^p_{kk',i,t}$}{Active power flow on line $k-k'$ under $S_{i,t}$ at $t$.}%
\Nom[P04]{$F^q_{kk',i,t}$}{Reactive power flow on line $k-k'$ under $S_{i,t}$ at $t$.}%
\Nom[P05]{$F^s_{kk'}$}{Apparent power capacity of line $k-k'$.}%
%\Nom[P06]{$G^p_{k,i,t}$}{Real power injection of DGs at bus $k$ (if any) under state $S_{i,t}$.}%
%\Nom[P07]{$G^q_{k,i,t}$}{Reactive power injection of DGs at bus $k$ (if any) under state $S_{i,t}$.}%
%\Nom[P08]{$\underline G^p_{k}, \overline G^p_{k}$}{Lower/upper limits of real power injection of DGs at bus $k$ (if any).}%
%\Nom[P09]{$\underline G^q_{k}, \overline G^q_{k}$}{Lower/upper limits of reactive power injection of DGs at bus $k$ (if any).}%
\Nom[P10]{$\Delta {L^p_{k,i,t}}, \Delta {L^q_{k,i,t}}$}{ Loss of active/reactive load of bus $k$ under state $S_{i,t}$.}%
\Nom[P11]{${L^p_{k,t}}$}{Active load of bus $k$ at $t$.}%
\Nom[P12]{${L^q_{k,t}}$}{Reactive load of bus $k$ at $t$.}%
\Nom[P13]{$M$}{Large positive number.}%
\Nom[P14]{$o_{kk',i,t}$}{Binary variable, the value is $1$ if bus $k'$ is the parent bus for bus $k$ under state $S_{i,t}$, otherwise $0$.}%
\Nom[P15]{$r_{kk'}/x_{kk'}$}{Resistance/reactance of line $k-k'$.}%
\Nom[P16]{$s_{c,t}$}{On-off state of component $c$ at $t$.}%
\Nom[P16]{$S_{i,t}, S_{j,t+1}$}{Markov state $i$ and $j$ at $t$ and $t+1$, respectively.}%
\Nom[P17]{$\Delta T$}{Duration of each time period.}%
\Nom[P18]{$T_c^f$}{Time period from normal state to failure state of component $c$.}%
\Nom[P19]{$\Delta T_c^r$}{Repair duration for component $c$.}%
\Nom[P20]{$U_{k,i,t}$}{Squared voltage magnitude of bus $k$.}%
\Nom[P21]{$\underline V_k, \overline V_k$}{Low/upper limits of voltage value of bus $k$.}%
\Nom[P22]{$v_t, v_{t+1}$}{values functions at $t$ and $t+1$.}%
\Nom[P23]{$\Pr$}{Transition probability.}%
\Nom[P24]{$\beta_{l,i,t}, \beta_{c,i,t}$}{Binary variables representing on-off states of line $l$ and $c$, respectively. $1$ denotes on state, and $0$ denotes off state.}%
\Nom[P25]{$\xi_t$}{Uncertainty of extreme event at $t$.}%
\Nom[P25]{$\eta_t, \eta_{t-1}$}{Penalty due to loss of load at $t$ and $t-1$ ($\$/MWh$).}%

\Nom[N01]{${\bf{b}}_{i}$}{Binary-coded matrix for post-decision states.}%
\Nom[N02]{$\mathbb{E}$}{Expected value.}%
\Nom[N03]{$m$}{Number of dispatchable lines.}%
\Nom[N04]{$n$}{Number of iterations.}%
\Nom[N05]{$N$}{Maximum number of iterations.}%
\Nom[N06]{$S_{i,t}^{a_t}$}{Post-decision after $S_{i,t}$ with action $a_t$.}%
\Nom[N06]{$S_{j',t-1}^{a_{t-1}}$}{Post-decision after $S_{j',t-1}$ with action $a_{t-1}$.}%
\Nom[N07]{$T$}{Number of decision periods.}%
\Nom[N08]{$v_t^{a_t}$}{Value function of post-decision state.}%
\Nom[N08]{$v_t^{n}$}{Value function of state at $n^{th}$ iteration at $t$.}%
\Nom[N09]{$\tilde v_t^{a_t,n}, \tilde v_{t-1}^{a_{t-1},n}$}{Approximated value function of post-decision state at $n^{th}$ iteration.}%
\Nom[N09]{$V_{i'}$}{Known values of post-decision states.}%
\Nom[N10]{$x_1,x_2$}{Binary variables.}%
\Nom[N11]{$y_{l,i',t}$}{Binary variable.}%
\Nom[N12]{$\epsilon$}{A coefficient.}%

%\Nom[P]{$S_{i,t}$}{State $i$ at decision epoch $t$.}%

%\Nom[N]{$\alpha_t$}{A vector representing an hyper-plane.}%
%\Nom[N]{$\Omega_{\Pi}^{Int}$}{Initial set of belief states.}%

\section{Introduction}
\IEEEPARstart{B}{ecause} distribution systems are directly connected to commercial and residential customers with radial topologies, any failures in distribution systems will lead to outages. Climate change increases the frequency and intensity of severe weather, which is a major cause of severe system failures. For example, weather events caused roughly $679$ power outages, each of which affected at least $50,000$ customers, between $2003$ and $2012$. Although transmission system outages did occur, a major portion of outages occurred along distribution systems \cite{AP_CosttoDistribution}. The severe consequences have required distribution systems to have resilience against these extreme weather events, and this has been identified by the United States Electric Power Research Institute (EPRI) \cite{AP_EPRI1} and the North American Electric Reliability Corporation (NERC) \cite{AP_NERC1}.

To enhance the system resilience, we can take actions in three stages with regard to weather events, i.e., prior to events, during events, and after events \cite{AP_ThreeStage1,  AP_ThreeStage2}. Prior to weather events, the historical data-based models \cite{AP_PredictModel1, AP_PredictModel2, AP_Duration1} are used to estimate outages that help the system operators to make preventive actions such as system maintenance \cite{AP_Maintenance1} and system hardening \cite{AP_Hardening1}. System hardening makes physical infrastructural changes to systems so that they are less susceptible to extreme events. For example, a coordinated hardening and distributed generator (DG) allocation strategy has been developed in \cite{AP_Hardening1}. Prior to events, preparing enough blackstart generators and emergency generators after potential failures is also a critical measure to improve the system resilience. To this end, a Generic Restoration Milestones (GRMs)-based algorithm is developed to assess blackstart capacities \cite{AP_Emergency2}, and a procurement plan with a minimal cost while guaranteeing sufficient blackstart capacities is proposed to provide enough blackstart resources at right locations \cite{AP_Emergency1}. To effectively isolate possible failures and connect blackstart/emergency generators to systems, a resilience-based model for switch placement in distribution systems is developed prior to events \cite{AP_SwitchAllocation1}. In addition to physical power systems, hardening  communication systems in charge of monitoring/controlling the physical power systems play an important role in enhancing the system resilience \cite{AP_WideControlRe1}. Even though many preventive actions are performed prior to events, it is impossible to avoid outages completely. When outages occur after events, it is necessary to recover outages as quickly as possible to improve the system resilience. A conventional power system restoration includes three stages, i.e., preparation, system restoration and load restoration \cite{AP_Restoration1, AP_Restoration2, AP_Restoration3}. Some algorithms such as expert systems \cite{AP_Restoration4} and heuristic approaches \cite{AP_Restoration5} are proposed to accelerate load recovery. However, there are unique characteristics associated with outages caused by weather-related events, leading to different restoration strategies such as microgrid-based restoration strategies \cite{AP_Restoration6} and  decentralized restoration schemes \cite{AP_Restoration7}.

The above studies mainly focus on strategies prior to events and after events. However, few studies investigate the strategies during the unfolding events. One difficulty in establishing strategies during the unfolding events is to map sequentially varying states caused by the unfolding events to optimal strategies. The commonly used scenario-based stochastic programming \cite{AP_StocBook1} and robust stochastic programming are not suitable for mapping sequentially varying states to optimal strategies. To address this difficulty, MDP can be employed to help make state-based decisions on a stochastic environment caused by weather events. Some applications of MDP in power systems have been investigated \cite{AP_MDPappl1, AP_MDPappl2}. For the resilience enhancement, \cite{AP_Chong1MDP} proposes sequentially proactive MDP-based  strategies to improve the transmission system resilience, and a linear scalarization method based on the state tree is used to solve the proposed model. However, distribution systems and transmission systems differ in topologies and allowable actions, the developed model and solution in \cite{AP_Chong1MDP} cannot be applied to distribution systems directly. It is necessary to develop state-based decision-making processes for distribution systems considering their own characteristics.   

This paper proposes MDP-based resilience enhancement for distribution systems. The contributions of this paper are three-fold: 1) The sequentially uncertain states, i.e., feeder line on/off states, on the trajectory of an unfolding event are represented as Markov states. Transition probabilities between Markov states are determined by component failure probabilities caused by the unfolding event. 2) A recursive optimization model for each Markov state is constructed to map states to optimal strategies. The allowable action for each state is system reconfiguration. 3) An approximate dynamic programming approach based on post-decision states and value function approximation is employed to solve the proposed model to deal with the curse of dimensionality caused by a mass of states and allowable actions.      

The remainder of this paper is organized as follows. Section II shows extreme events' impacts on system states. Section III presents the mathematical formulation, and section IV introduces the solution method. The case studies are demonstrated in Section V, and the work is concluded in Section VI. 

\section{Modeling Influences of Events on Distribution Systems}
This section first introduces Markov states on the trajectories of extreme events, and then presents transition probabilities between different Markov states under allowable actions.
\subsection{Markov states on the event's trajectories}
Usually, the impacts of a weather-related event on a distribution system are sequential due to the sequential trajectory, indicating that the components with different locations in the system may be in failure in different time periods. This results in the sequential changes of system states such as on-off states of distribution lines. A system state including on-off states of distribution lines on the trajectory is represented as a Markov state in this paper. Define $\mathcal{F}_t$ as the set of all possible failure components due to the unfolding event at $t$, and $\tilde{\mathcal{F}}_t$ as the set of the actual failure components at $t$, and we have $\tilde{\mathcal{F}}_t \subseteq \mathcal{F}_t$. For example, $\mathcal{F}_{t_2} = \{ b_{1-2}, b_{1-4} \}$ and $\tilde {\mathcal{F}}_{t_2} = \{ b_{1-4} \}$. The failure scenarios of the lines $b_{1-2}$ and $b_{1-4}$ are uncertain before the time period $t_2$, and the actual state can only be observed at $t_2$ and the actual failure on $b_{1-4}$ occurs. Considering the sequential characteristic of the extreme event, the Markov state $S_{i,t}$ at $t$ can be represented as follows.
\begin{align}
& {S_{i,t}} = \bigcup\limits_{\tau = 1}^t {\left( \tilde{\mathcal{F}}_{\tau}- \tilde{\mathcal{R}}_{\tau} \right) } \label{StateUnion}
\end{align}
where (\ref{StateUnion}) shows that the Markov state at $t$ is represented as failure components except repaired components from the initial time period to the time period $t$. $\tilde{\mathcal{R}}_{\tau}$ is the set of components repaired at time $\tau$.  

\subsection{Transition probability between Markov states}
On the trajectory of the event, the current Markov state at $t$ has a probability of reaching to each future Markov state at $t+1$, and the probability is called as a transition probability, which is determined by component failure rates caused by the event and the allowable system reconfiguration. The transition probability can be expressed as follows.
\begin{align}
& \Pr ({S_{j,t + 1}}|{S_{i,t}},{a_t},\xi_t ) = \prod\nolimits_{c \in \mathcal{F}_{t+1} } {\Pr ({s_{c,t + 1}}|{s_{c,t}},{a_t},\xi_t )} 
\label{TransPro}
\end{align}    
where $\Pr ({S_{j,t + 1}}|{S_{i,t}},{a_t},\xi_t )$ means the probability from the state ${S_{i,t}}$ to the state ${S_{j,t + 1}}$ under the action ${a_t}$ with the uncertainty $\xi_t$, and ${\Pr ({s_{c,t + 1}}|{s_{c,t}},{a_t},\xi_t )}$ represents the probability from the $c$th component's on-off state $s_{c,t}$ to the on-off state $s_{c,t+1}$ under the action ${a_t}$ with the uncertainty $\xi_t$. Three causes can change the on-off states, and they are listed as follows.
\begin{itemize}
	\item  Component failure caused by extreme events: Since the on-off states are uncertain in the next time period due to extreme events, it is necessary to calculate the probability of each scenario. At present, there are many existing studies on components' failure probabilities caused by extreme events  \cite{Failure111}.
	\item  System reconfiguration: After system reconfiguration, the line state, i.e., on state or off state, is determined. This indicates that the corresponding probability from $s_{k,t}$ to $s_{k,t+1}$ is $0$ or $1$.  
	\item  Repair: Before the failure components are repaired, the state is off. After repaired, the state is on. This shows that the corresponding probability from $s_{k,t}$ to $s_{k,t+1}$ is $0$ or $1$.
\end{itemize} 
\begin{figure}[!h]
	\centering
	\includegraphics[width=5cm]{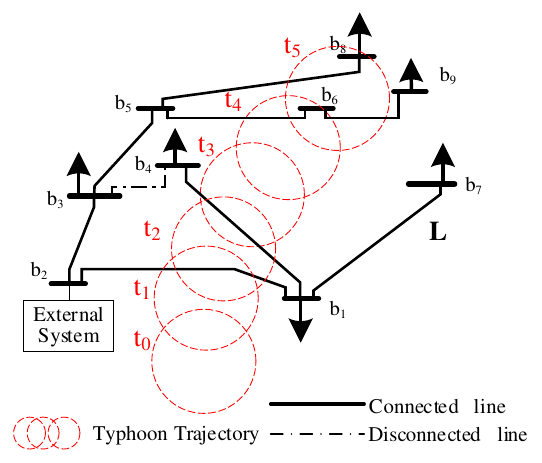}
	\caption{An example of a distribution system under an unfolding event}
	\label{StateDemonstrate_fig}
\end{figure}

\section{Optimization Model based on Markov Decision Processes}
This section first introduces a recursive model to map each Markov state to its optimal strategy, and then lists the operational constraints for distribution systems.
\subsection{Markov Decision Processes-based Recursive Model}
Transition probabilities in (\ref{TransPro}) show that the current actions associated with uncertainties caused by extreme events impact the current states and future states. Since different states result in different operational costs, it is necessary to make decisions based on not only current states but also future states impacted by transition probabilities. Take the scenario in Fig. \ref{StateDemonstrate_fig} as an example. Since the line $b_{1-2}$ will be impacted by the typhoon at $t_1$, the line $b_{1-2}$ may be in failure at $t_1$, resulting in the outages of the lines $b_{1-4}$ and $b_{1-7}$. If we can disconnect the line $b_{1-2}$ and connect the line $b_{3-4}$ before the line $b_{1-2}$ is impacted by the typhoon, we can avoid the outages of the lines $b_{1-4}$ and $b_{1-7}$ at $t_1$. For this case, we only consider one time period ahead and there is only one line impacted by the typhoon at $t_1$. If all time periods over the unfolding typhoon and numerous components impacted in each time period are considered, we need to develop a model to help make decisions to ensure the minimum operational cost. Considering the sequential time periods and numerous components on the trajectory, we establish a recursive model with the current cost and the expected future cost listed as follows.  
\begin{align}
&\begin{array}{l}
{v_t}({S_{i,t}}) = \\[5pt]
\mathop {\min }\limits_{{a_t} \in \mathcal{A}} \left( \begin{array}{l}
{C_t}({S_{i,t}},{a_t}) + \\[5pt]
\sum\limits_{{S_{j,t + 1}} \in \mathcal{S}} {\Pr ({S_{j,t + 1}}|{S_{i,t}},{a_t},{\xi _t}) \cdot {v_{t + 1}}({S_{j,t + 1}})} 
\end{array} \right)
\end{array}
\label{RecursiveModel}
\end{align}
where ${v_t}({S_{i,t}})$ and ${v_{t+1}}({S_{j,t+1}})$ are the value functions of the states $S_{i,t}$ and $S_{j,t+1}$ at $t$ and $t+1$, respectively. The second term on the right side of (\ref{RecursiveModel}) shows the expected future cost. ${C_t}({S_{i,t}},{a_t})$ is the current cost caused by the action $a_t$ under the state $S_{i,t}$ at $t$, and this cost in this study is defined as the sum of the cost of loss of load and the operational cost of controllable lines. It is expressed as follows. 
\begin{align}
&{C_t}({S_{i,t}},{a_t}) = \sum\limits_{k \in \mathcal{B}} {({\eta _t} \cdot \Delta {L^p_{b,i,t}} \cdot \Delta T)}  + \sum\limits_{l \in \mathcal{L}_{i,t}^d} {({\beta _{l,i,t}} \cdot {C_l})} 
\label{CostExpre}
\end{align}
where the first term on the right side of (\ref{CostExpre}) is the cost of loss of load, and the second term is the operational cost of controllable lines.  

\subsection{Operational Constraints} 
For each time period during the unfolding event, the operational constraints, i.e., radial topologies, power balance, power flow, voltage limits, and line capacity, should be satisfied.
\subsubsection{Radiality constraint} Different from transmission systems, distribution systems should operate in radial topologies. When performing system reconfiguration under the state $S_{i,t}$, the spanning tree constraints are used to guarantee the network radiality.  
\begin{subequations}
	\label{RadialityC}
	\begin{align}
	& {o_{kk',i,t}} + {o_{k'k,i,t}} = {\beta _{l,i,t}}\quad l \in {\mathcal{L}_{i,t}},t \in \mathcal{T} \label{RadialityC1}\\[2pt]
	& \sum\limits_{k' \in {\mathcal{N}_{k,i,t}}} {{o_{kk',i,t}}}  = 1\quad k \in {{ \mathcal{\bar B}}_{i,t}},t \in \mathcal{T} \label{RadialityC2}\\[2pt]
	& {o_{kk',i,t}} = 0\quad k \in \mathcal{\tilde B},k' \in {\mathcal{N}_{k,i,t}},t \in \mathcal{T} \label{RadialityC3}
	\end{align}
\end{subequations}
where (\ref{RadialityC1}) and (\ref{RadialityC2}) constrain that the two terminals of a connected line only have one parent bus. In practice, islanded buses, to which power cannot be supplied by the grid, maybe exist due to component failures caused by extreme events, and these islanded buses are not included in the spanning tree constraint. (\ref{RadialityC3}) indicates that the substation bus (i.e., the bus connected to the external system) has no parent buses.

In practice, it is possible that only some lines can be dispatched under the state $S_{i,t}$. In this case, we can add an constraint with regard to non-dispatched lines.
\begin{align}
& {\beta _{l,i,t}} = 1\quad l \in {\mathcal{L}^{nd}_{i,t}}, t \in \mathcal{T}
\label{NonDisline}
\end{align}

\subsubsection{Repair constraint} When there are components in failure under the state $S_{i,t}$, these failure components and the resulting islanded components cannot participate in system dispatch before they are repaired.
\begin{align}
& {\beta _{c,i,t}} = 0\quad c \in {\mathcal{C}^{f}_{i,t}},T_c^{f} \le t \le T_c^{f} + \Delta T_c^r
\label{ReCon1}
\end{align}
where (\ref{ReCon1}) means that the state of the failure component $c$ is set to $0$ during the repair time periods.

\subsubsection{Power flow constraint}
The power flow of each line has relations to bus voltages of the two terminal buses of each line, and can be expressed as follows.
\begin{subequations}
	\label{PowerFlow}
	\begin{align}
	&\begin{array}{l}
	{U_{k,i,t}} - {U_{k',i,t}} \le (1 - {\beta _{l,i,t}}) \cdot M + \\[4pt]
	\quad  2({r_{kk'}} \cdot F_{kk',i,t}^p + {x_{kk'}} \cdot F_{kk',i,t}^q)
	\quad l \in {\mathcal{L}_{i,t}},t \in \mathcal{T}
	\end{array} \label{PowerFlow1}\\[3pt]
	&\begin{array}{l}
	{U_{k,i,t}} - {U_{k',i,t}} \ge ({\beta _{l,i,t}} -1 ) \cdot M + \\[4pt]
	\quad 2({r_{kk'}} \cdot F_{kk',i,t}^p + {x_{kk'}} \cdot F_{kk',i,t}^q)
	\quad l \in {\mathcal{L}_{i,t}},t \in \mathcal{T}
	\end{array} \label{PowerFlow2}
	\end{align}
\end{subequations}
where (\ref{PowerFlow1})-(\ref{PowerFlow2}) are derived from the DsitFlow model \cite{DistFlow1}. The quadratic terms in the accurate power flow model are ignored \cite{DistFlow2}. The big $M$ is a disjunctive parameter. With a sufficiently large $M$, (\ref{PowerFlow1})-(\ref{PowerFlow2}) are redundant when distribution lines are disconnected or outages.  Non-islanded buses are included in these constraints.

\subsubsection{Power balance constraint}
When reaching to the state $S_{i,t}$ at $t$, the out-flow/in-flow power of each non-islanded bus in the system should be equal. The constraint can be expressed as follows.  
\begin{subequations}
	\label{PowerBan}
	\begin{align}
	&\begin{array}{l} \left( {{L^p_{k,t}} - \Delta {L^p_{k,i,t}}} \right) + \sum\limits_{k' \in {\mathcal{N}_{k,i,t}}} {{F^p_{kk',i,t}}}  = 0 \quad \\
	\quad\quad\quad\quad\quad\quad\quad\quad\quad\quad\quad\quad\quad\quad\quad k \in {{ \mathcal{\bar B}}_{i,t}},t \in \mathcal{T}  
	\end{array} \label{PowerBan1}\\[2pt]
	&\begin{array}{l} \left( {{L^q_{k,t}} - \Delta {L^q_{k,i,t}}} \right) + \sum\limits_{k' \in {\mathcal{N}_{k,i,t}}} {{F^q_{kk',i,t}}}  = 0 \quad \\
	\quad\quad\quad\quad\quad\quad\quad\quad\quad\quad\quad\quad\quad\quad\quad k \in {{ \mathcal{\bar B}}_{i,t}},t \in \mathcal{T}  
	\end{array} \label{PowerBan2}
	\end{align}
\end{subequations}  
where (\ref{PowerBan1}) and (\ref{PowerBan2}) represent real power balance and reactive power balance, respectively. Only non-islanded buses are included in the constraint. The load connected to the islanded buses in the system is directly considered as loss of load in (\ref{CostExpre}).

\subsubsection{Line capacity constraint}
The power through each line should be within the limit for the state $S_{i,t}$ with the action $a_t$. The constraint can be expressed as follows.
\begin{align}
& {(F_{kk',i,t}^p)^2} + {(F_{kk',i,t}^q)^2} \le {\beta _{l,i,t}} \cdot {(F_{kk'}^s)^2} \quad l \in {\mathcal{L}_{i,t}},t \in \mathcal{T}
\label{LineCap1}
\end{align}
where (\ref{LineCap1}) is a nonlinear constraint, resulting in computational intractability. To facilitate the model solution, the constraint (\ref{LineCap1}) is relaxed to a group of linear constraints \cite{LinePower1}, and are rewritten as follows.  
\begin{subequations}
	\label{LineCap_linear}
	\begin{align}
	& \begin{array}{l}
	 - {\beta _{l,i,t}} \cdot F_{kk'}^s \le F_{kk',i,t}^p \le {\beta _{l,i,t}} \cdot F_{kk'}^s \\[3pt]
	 \quad\quad\quad\quad\quad\quad\quad\quad\quad\quad\quad\quad\quad\quad\quad  l \in {\mathcal{L}_{i,t}},t \in \mathcal{T}
	 \end{array}  \label{LineCap_linear1}\\[2pt]
	& \begin{array}{l}
	- {\beta _{l,i,t}} \cdot F_{kk'}^s \le F_{kk',i,t}^q \le {\beta _{l,i,t}} \cdot F_{kk'}^s \\[3pt]
	\quad\quad\quad\quad\quad\quad\quad\quad\quad\quad\quad\quad\quad\quad\quad  l \in {\mathcal{L}_{i,t}},t \in \mathcal{T}
	\end{array}  \label{LineCap_linear2}\\[2pt]
	& \begin{array}{l}
	- \sqrt 2 {\beta _{l,i,t}} \cdot F_{kk'}^s \le F_{kk',i,t}^p + F_{kk',i,t}^q \le \sqrt 2 {\beta _{l,i,t}} \cdot F_{kk'}^s \\[3pt]
	\quad\quad\quad\quad\quad\quad\quad\quad\quad\quad\quad\quad\quad\quad\quad\quad  l \in {\mathcal{L}_{i,t}},t \in \mathcal{T}
	\end{array}  \label{LineCap_linear3}\\[2pt]
	& \begin{array}{l}
	- \sqrt 2 {\beta _{l,i,t}} \cdot F_{kk'}^s \le F_{kk',i,t}^p + F_{kk',i,t}^q \le \sqrt 2 {\beta _{l,i,t}} \cdot F_{kk'}^s \\[3pt]
	\quad\quad\quad\quad\quad\quad\quad\quad\quad\quad\quad\quad\quad\quad\quad\quad  l \in {\mathcal{L}_{i,t}},t \in \mathcal{T}
	\end{array}  \label{LineCap_linear4}
	\end{align}
\end{subequations}  
   
\subsubsection{Voltage constraint} The voltage limits under the state $S_{i,t}$ with the action $a_{t}$ should be satisfied.
\begin{align}
& \underline V_k^2 \le {U_{k,i,t}} \le \overline V_k^2 \quad k \in {\mathcal{\bar B}_{i,t}},t \in \mathcal{T}
\label{Vlimits}
\end{align}

\subsection{MDP-based Optimization Model}
The constructed MDP-based optimization model can be represented as follows.
\begin{align}
& \begin{array}{l}
{\rm{Obj. \quad\quad \;\;\;\;\;        (\ref{RecursiveModel})}}\\
{\rm{s.t.\quad \quad \quad \quad \,                        (\ref{CostExpre})-(\ref{Vlimits})}}
\end{array} 
\label{Equi_Model1}
\end{align} 

This is a recursive model for each Markov state which is computationally intractable. The next section will introduce the solution method. 

%CostExpre

\section{Model Solution based on Approximate Dynamic Programming}
This section first introduces the challenges of solving the proposed model, and then presents the basic idea of ADP in solving the MDP-based model, and finally shows how to solve the proposed model by using the ADP approach.
\subsection{Challenges of model solution}
When using the conventional stochastic programming to deal with a sequential decision-making problem with uncertainties, one commonly used approach is to generate some scenarios to represent the uncertainties. Based on these generated scenarios, we optimize a model with an expected objective, and then we can obtain the optimal strategy. The conventional stochastic programming cannot be used to deal with the proposed model because the influences of actions on scenario transitions are not included and the real scenario may not be included in the generated scenarios. 

\begin{figure*}
	\centering
	\includegraphics[width=16cm]{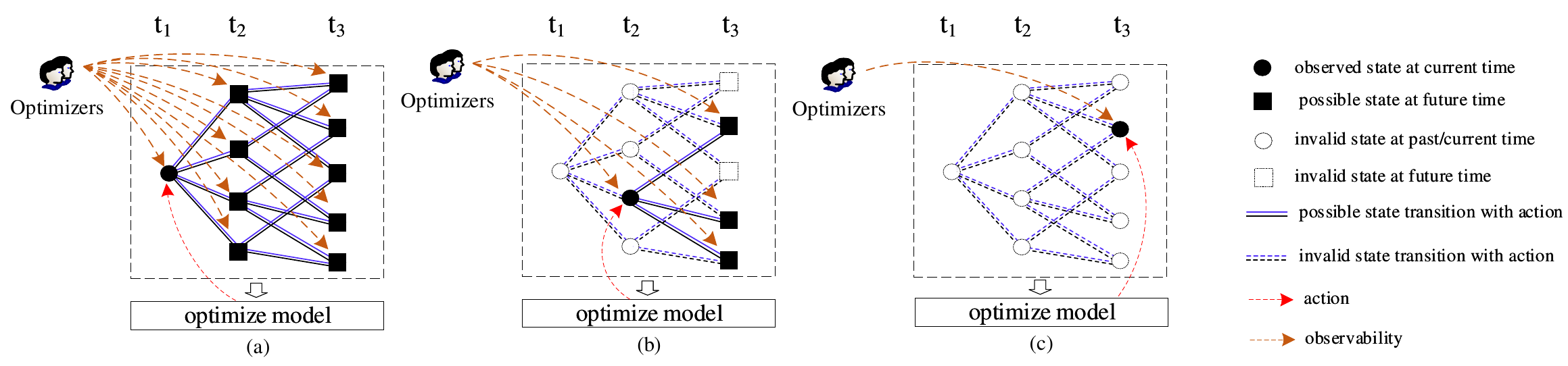}
	\caption{Markov state-based decision at (a) $t_1$, (b) $t_2$, and (c) $t_3$.}
	\label{ChallengeMDP_fig}
\end{figure*}

The decision processes used in this paper can be illustrated by using the case in Fig. \ref{ChallengeMDP_fig} (a), (b), and (c), which represent the decision making processes at $t_1$, $t_2$, and $t_3$, respectively. At $t_1$, the state transition tree under the impacts of actions is constructed, as shown in Fig. \ref{ChallengeMDP_fig} (a), and the optimal action can be obtained by optimizing the recursive model (\ref{RecursiveModel}). After performing the optimal action, the state reaches to a new state at $t_2$ under uncertainty, as shown in Fig. \ref{ChallengeMDP_fig} (b). For this new state, the state transition tree under the impacts of actions needs to be updated again because some state transitions may be invalid. With the new state transition tree, the optimal action for the new state at $t_2$ can be obtained by optimizing the recursive model (\ref{RecursiveModel}) again. The optimal action for the new state at $t_3$ can be obtained with the similar processes. According to the decision processes, constructing the state transition tree under actions in consideration of uncertainties is one critical step to solve the MPD-based model. However, it is a difficult task to construct the state transition tree of the proposed model in consideration of various actions and the resulting complicated state transitions. In addition, a large-scale problem with numerous states is possibly intractable due to ``curse of dimensionality''. There are three curses of dimensionality: (i) the state space $\mathcal{S}$ may be too huge to calculate the value function ${v_t}({S_{i,t}})$ for each state within acceptable time, (ii) the decision space $\mathcal{A}$ is too large to obtain the optimal action for each state, (iii) the outcome space may be too large to calculate the expectation of future cost.         
\subsection{Approximate dynamic programming}
Approximate dynamic programming is a modeling framework offering some techniques for dealing with the curses of dimensionality in multi-period, large, and stochastic MDP-based models. There are two critical techniques used by ADP: 1) post-decision states are constructed to deal with the large outcome space, 2) a forward dynamic algorithm based on sample paths is used to solve the recursive model by stepping forward in time, and repeat this procedure for enough iterations. 
\subsubsection{Post-decision states} The post-decision state, defined as ${S_{i,t}^{a_t}}$, is a state immediately after the action $a_t$ but before the arrival of a new state in consideration of uncertainties. To apply the ADP approach, a more generic form of the proposed model in (\ref{RecursiveModel}) is an expectational form listed as follows.  
\begin{align}
& {v_t}({S_{i,t}}) = \mathop {\min }\limits_{{a_t} \in \mathcal{A}} \left( {{C_t}({S_{i,t}},{a_t}) + \mathbb{E} \{ {v_{t + 1}}({S_{j,t + 1}}|{S_{i,t}},{a_t},{\xi _t})\} } \right)
\label{GenericModel1}
\end{align}
where (\ref{GenericModel1}) can be rewritten as (\ref{GenericModel2}) with the post-decision state ${S_{i,t}^{a_t}}$.  
\begin{align}
& {v_t}({S_{i,t}}) = \mathop {\min }\limits_{{a_t} \in \mathcal{A}} \left( {{C_t}({S_{i,t}},{a_t}) + \mathbb{E} \{ {v_{t + 1}}({S_{j,t + 1}}|{S_{i,t}^{a_t}},{\xi _t})\} } \right)
\label{GenericModel2}
\end{align}

Define ${\mathbb{E} \{ {v_{t + 1}}({S_{j,t + 1}}|{S_{i,t}^{a_{t}}},{\xi _t})\}}$ by $v_t^{a_{t}}({S_{i,t}^{a_{t}}})$, we have the following optimality equations.
\begin{subequations}
	\label{New_GenericModel}
	\begin{align}
	& {v_t}({S_{i,t}}) = \mathop {\min }\limits_{{a_t} \in \mathcal{A}} \left( {{C_t}({S_{i,t}},{a_t}) + v_t^{a_t}(S_{i,t}^{a_t})} \right)
	\label{New_GenericModel1}\\[2pt]
	& v_{t-1}^{a_{t-1}}({S_{j',t-1}^{a_{t-1}}}) = {\mathbb{E} \{ {v_{t}}({S_{i,t}}|{S_{j',t-1}^{a_{t-1}}},{\xi _{t-1}})\}}
	\label{New_GenericModel2}
	\end{align}
\end{subequations}  

Substituting (\ref{New_GenericModel1}) into (\ref{New_GenericModel2}) results in the optimality equations of the post-decision states as follows.
\begin{align}
& v_{t - 1}^{{a_{t - 1}}}(S_{j',t - 1}^{{a_{t - 1}}}) = \mathbb{E} \{ \mathop {\min }\limits_{{a_t} \in \mathcal{A}} ({C_t}({S_{i,t}},{a_t}) + v_t^a(S_{i,t}^a|S_{j',t - 1}^{{a_{t - 1}}},{\xi _{t - 1}}))\} 
\label{PostDecisiont_1}
\end{align}  
where (\ref{PostDecisiont_1}) can be rewritten as the form at $t$ as follows. 
\begin{align}
& \begin{array}{l} 
v_{t}^{{a_{t}}}(S_{i,t}^{{a_{t}}}) = \\[2pt]
\quad \mathbb{E} \{ \mathop {\min }\limits_{{a_{t+1}} \in \mathcal{A}} ({C_{t+1}}({S_{j,t+1}},{a_{t+1}}) + v_{t+1}^{a_{t+1}}(S_{j,t+1}^{a_{t+1}}|S_{i,t}^{{a_{t}}},{\xi _{t}}))\} 
\end{array}
\label{PostDecisiont}
\end{align}
where (\ref{PostDecisiont}) shows the value of the post-decision state.

When $v_{t}^{{a_{t}}}(S_{i,t}^{{a_{t}}})$ in (\ref{New_GenericModel1}) is known, it would be easy to solve the optimization model (\ref{New_GenericModel1}). Based on this idea, ADP is to use the deterministic optimization model (\ref{New_GenericModel1}) with an initial estimation of $\tilde v_t^{a_t}(S_{i,t}^{a_t})$ of $v_t^{a_t}(S_{i,t}^{a_t})$ to make decisions for each state, and then employ the resulting observations to update an estimation $\tilde v_t^{a_t}(S_{i,t}^{a_t})$ thereby approximating the expected value in (\ref{PostDecisiont}).      
   
\subsubsection{Forward Dynamic algorithm} 
For the forward dynamic algorithm, the recursive model is solved only for one state in each time period, by using the estimation of the post-decision state and performing iterations to update the estimations of the post-decision states on the sample paths. To deal with iterations, we add a superscript $n$ and $n-1$ to the value functions, and (\ref{New_GenericModel1}) can be expressed as follows.   
\begin{align}
& {v^n_t}({S_{i,t}}) = \mathop {\min }\limits_{{a_t} \in \mathcal{A}} \left( {{C_t}({S_{i,t}},{a_t}) + \tilde v_t^{a_t,n-1}(S_{i,t}^{a_t})} \right)
\label{Iteration_model1}
\end{align}
where the decision that minimizes (\ref{Iteration_model1}) at $n$th iteration is shown as follows.
\begin{align}
& a_t = \mathop {\arg \min }\limits_{{a_t} \in \mathcal{A}} \left( {{C_t}({S_{i,t}},{a_t}) + \tilde v_t^{a_t,n-1}(S_{i,t}^{a_t})} \right)
\label{Iteration_model2}
\end{align} 
%where $\tilde v_t^{a_t,n}(S_{i,t}^{a_t})$ can be updated by 
%\begin{align}
%& \tilde v_t^{a_t,n}(S_{i,t}^{a_t}) = (1-\epsilon) \cdot \tilde v_t^{a_t,n-1}(S_{i,t}^{a_t}) + \epsilon \cdot {v^n_{t+1}}({S_{j,t+1}}) 
%\label{Update_model}
%\end{align} 
where $\tilde v_t^{a_t,n-1}(S_{i,t}^{a_t})$ can be updated by 
\begin{align}
& \tilde v_t^{a_t,n-1}(S_{i,t}^{a_t}) = (1-\epsilon) \cdot \tilde v_t^{a_t,n-2}(S_{i,t}^{a_t}) + \epsilon \cdot {v^{n-1}_{t+1}}({S_{j,t+1}}) 
\label{Update_model}
\end{align} 
where the first term on the right side of (\ref{Update_model}) represent the estimate of the post-decision state $S_{i,t}^{a_t}$ at the $(n-2)$th iteration, and the second term represent the value of the resulting observations from the post-decision state $S_{i,t}^{a_t}$ at the $n$th iteration. 

\subsection{Reformulation of the proposed model}
Based on (\ref{Iteration_model1}), we just need to solve a deterministic model. In the model, the term $C_t(S_{i,t},a_t)$ is an explicit objective (\ref{CostExpre}) with regard to variables associated with constraints (\ref{RadialityC})-(\ref{Vlimits}), however, the term $\tilde v_t^{a_t,n}(S_{i,t}^{a_t})$ is just a value with regard to the post-decision state $S_{i,t}^{a_t}$ but has no relations to the variables and actions. In this case, it is not possible to optimize the model (\ref{Iteration_model1}). Therefore, it is necessary to relate $\tilde v_t^{a_t,n}(S_{i,t}^{a_t})$ to the variables and actions. 

%LineCap1

The Markov state in the study is determined by the on-off states of distributed lines. Since a failure is an observed event and the repair is an activity with continuous time period, whether a failure component is repaired or not at the current period is known. Furthermore, the post-decision states at the current period are defined as states before arrival of uncertainties in the next time period. Therefore, system reconfiguration is the cause of changing the current state $S_{i,t}$ to $S_{i,t}^{a_t}$. It is assumed that there are two reconfigurable lines and the corresponding binaries are $x_1$ and $x_2$. We will have four post-decision states listed in Table \ref{table1}. In this case, the second term $\tilde v_t^{a_t,n}(S_{i,t}^{a_t})$ in (\ref{Iteration_model1}) can be expressed as $(1-x_1)(1-x_2)V_1+(1-x_1)x_2V_2+x_1(1-x_2)V_3+x_1x_2V_4$ which relates the values of post-decision states to decision variables.   
\begin{table}[H]
	\centering
	\renewcommand{\arraystretch}{1.3}
	\caption{My caption}
	\label{table1}
	\begin{tabular}{p{2.6cm}<{\centering}p{0.5cm}<{\centering}p{0.5cm}<{\centering}p{2.6cm}<{\centering}}
		\hline
		Post-decision states    & $x_1$ & $x_2$ & Estimated Values   \\ \hline
		$S_{1,t}^{a_t}$                   & $0$   & $0$   & $\tilde v_t^{a_t,n}(S_{1,t}^{a_t}) \buildrel \Delta \over = V_1  $  \\ 
		$S_{2,t}^{a_t}$                   & $0$   & $1$   & $\tilde v_t^{a_t,n}(S_{2,t}^{a_t}) \buildrel \Delta \over = V_2  $  \\ 
		$S_{3,t}^{a_t}$                   & $1$   & $0$   & $\tilde v_t^{a_t,n}(S_{3,t}^{a_t}) \buildrel \Delta \over = V_3  $ \\ 
		$S_{4,t}^{a_t}$                   & $1$   & $1$   & $\tilde v_t^{a_t,n}(S_{4,t}^{a_t}) \buildrel \Delta \over = V_4  $ \\ \hline
	\end{tabular}
\end{table} 

Based on this technique, we can rewrite the second term $\tilde v_t^{a_t,n}(S_{i,t}^{a_t})$ in (\ref{Iteration_model1}) in a generic form. For the state $S_{i,t}$, there are $m$ reconfigurable lines, resulting in $2^m$ post-decision states associated with the corresponding estimates (represented as $V_1, V_2, \cdots, V_{2^m}$) at the $n$th iteration. The $2^m$ post-decision states are binary-coded with a $2^m \times m$ matrix ${\bf{b}}_i$, in which the $(i',l)^{th}$ entry denotes the on-off state of the line $l$ under the ${i'}^{th}$ post-decision states. The generic form of the term $\tilde v_t^{a_t,n}(S_{i,t}^{a_t})$ at the $n^{th}$ iteration is listed as follows. 
\begin{align}
& \sum\limits_{i' \in {\mathcal{S}_{i}^{post}}} {\left\{ {\prod\limits_{l \in {\mathcal{L}_{i,t}^{d}}} {(1 - {\beta _{l,i',t}} - {{\bf{b}}_i}(i',l))(1 - 2{{\bf{b}}_i}(i',l)){V_{i'}}} } \right\}} 
\label{GenericM_postD}
\end{align}      
where ${\beta _{l,i',t}}$ is a binary representing on-off states of the line $l$ under post-decision states. ${{\bf{b}}_i}(i',l)$ is a known value with 0 or 1, making (\ref{GenericM_postD}) a sum of multilinear functions. The optimization model (\ref{Iteration_model1}) in forward dynamic algorithm can be rewritten as follows.
\begin{align}
& \begin{array}{l}
{\rm{min\quad \quad \quad \;\;\;        (\ref{CostExpre})+(\ref{GenericM_postD})}}\\
{\rm{s.t.\quad \quad \quad \quad \,                        (\ref{RadialityC}),(\ref{NonDisline}),(\ref{ReCon1}),(\ref{PowerFlow}), (\ref{PowerBan}), (\ref{LineCap1}),(\ref{LineCap_linear}),(\ref{Vlimits})}}
\end{array} 
\label{Equi_Model2}
\end{align}

For multilinear functions in (\ref{GenericM_postD}), McCormick proposed a recursive procedure in which additional variables and constraints are added to obtain a formulation of the problem having only bilinear equations, which can be represented by four binary inequations. In (\ref{GenericM_postD}), the multilinear function with most variables is $\beta _{1,i',t}\beta _{2,i',t} \cdots \beta _{m,i',t}$, which can be represented by additional variables  
\begin{align}
\begin{array}{l}
{y_{2,i',t}} = {\beta _{1,i',t}}{\beta _{2,i',t}}\\
{y_{3,i',t}} = {y_{2,i',t}}{\beta _{3,i',t}}\\
\cdots \\
{y_{m,i',t}} = {y_{m - 1,i',t}}{\beta _{m,i',t}}
\end{array}
\label{NewVariables}
\end{align}
and additional constraints 
\begin{align}
\begin{array}{l}
{y_{2,i',t}} \ge {\beta _{2,i',t}} + {\beta _{1,i',t}} - 1\\
{y_{2,i',t}} \le {\beta _{1,i',t}}\\
{y_{l,i',t}} \ge 0\quad (l = 2, \cdots ,m)\\
{y_{l,i',t}} \le {\beta _{l,i',t}}\quad (l = 2, \cdots ,m)\\
{y_{l,i',t}} \ge {\beta _{l,i',t}} + {y_{l - 1,i',t}} - 1\quad (l = 3, \cdots ,m)\\
{y_{l,i',t}} \le {y_{l - 1,i',t}}\quad (l = 3, \cdots ,m)
\end{array}
\label{NewConstraints}
\end{align}
where (\ref{NewConstraints}) is an exact reformulation of $\beta _{1,i',t}\beta _{2,i',t} \cdots \beta _{m,i',t}$ since $\beta _{l,i',t}, l={1,2, \cdots, m}$ are binary variables. Based on the additional variables and the additional constraints, the optimization model is a mixed integer linear programming, which can be solved by many solvers such as CPLEX and GUROBI.   

\subsection{Solution procedure}
Based on Section IV.B and Section IV.C, we can solve the proposed MDP-based model by means of the
iteration-based ADP approach, and the CPLEX solver is used to solve the MILP model at each iteration. The detailed procedure is listed in {\bf{Algorithm 1}}. 
\begin{algorithm}
	\caption{ADP algorithm}
	{\small{
	\begin{algorithmic}[1] %sssssssssssssssss
		\State {Step 1. Initialization}
		\State {\quad\; Step 1.1. Set the iteration counter $n=1$ and the maximum number of iterations $N$.}
		\State {\quad\; Step 1.2. Set the initial approximation for each state.}
		
		\vspace{1mm}
		\State {Step 2. Do for $t=1, \cdots, T$}\vspace{1mm}
		
		\State {\quad\; Step 2.1. Solve (\ref{Iteration_model1}) and (\ref{Iteration_model2}) to get $v_t^n(S_{i,t})$ and $a_t$ at the iteration $n$. To solve (\ref{Iteration_model1}), (\ref{Equi_Model2}) and (\ref{NewConstraints}) will be used.}\vspace{1mm}
		\State {\quad\; Step 2.2. Update the approximation $\tilde v_t^{a_t,n}(S_{i,t}^{a_t})$ for the post-decision $S_{i,t}^{a_t}$  with (\ref{Update_model}).}\vspace{1mm}
		\State {\quad\; Step 2.3. Obtain the post-decision $S_{i,t}^{a_t}$ from the state $S_{i,t}$ under $a_t$.}\vspace{1mm}
		\State {\quad\; Step 2.4. According to uncertainties of the extreme event, generate a new state $S_{j,t+1}$ at $t+1$ from the post-decision state $S_{i,t}^{a_t}$ at $t$.}
		
		\vspace{1mm}
		\State {Step 3. Increment n. If $n \le N$ go to Step 1.}
		
		\vspace{1mm}
		\State {Step 4. Return the value function $\tilde v_t^{a_t,N}(S_{i,t}^{a_t})$ and the corresponding action $a_t$ for the state $S_{i,t}$.}

	\end{algorithmic} }}
\end{algorithm}

\section{Case Studies}
In this section, two test systems are used to verify the proposed model and the algorithm. The first system is the IEEE 33-bus system, and the second system is the IEEE 123-bus system. The cases are tested in MATLAB 2017a using  the CPLEX 12.6 solver on computers with 3.1 GHz i5 processors and 8 GB RAMS.

\subsection{IEEE 33-bus system}
\subsubsection{Data description}
\begin{figure}[!h]
	\centering
	\includegraphics[width=6cm]{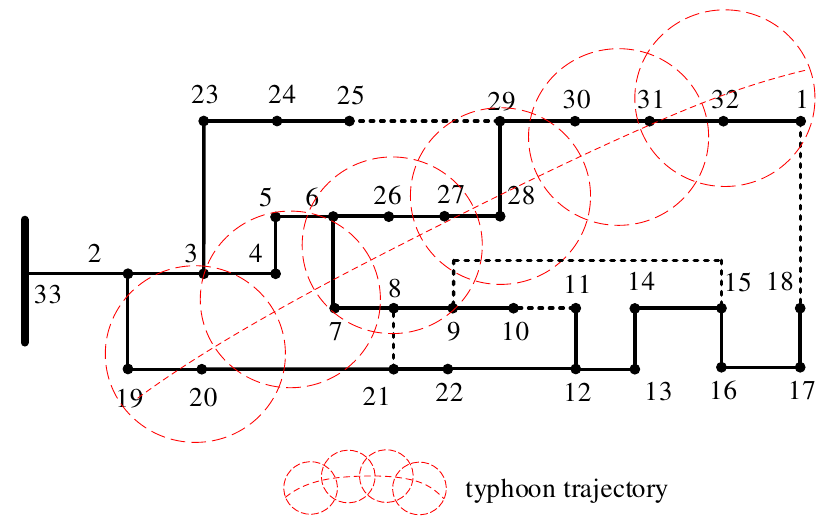}
	\caption{Topology of IEEE 33-bus system.}
	\label{IEEE33_fig}
\end{figure}

Fig. \ref{IEEE33_fig} shows the topology of the IEEE 33-bus system. The typhnoon trajectory is also shown in Fig. \ref{IEEE33_fig}. For the original topology, the lines 8-21, 12-22, 1-18, 9-15, and 25-29 are disconnected to ensure the radial topology. It is assumed that the lines 10-11, 12-13, 25-29, 1-18, 14-15, 12-22, 8-21, and 9-15 are dispatchable and the other lines are non-dispatchable. %The operational costs for the lines 10-11, 12-13, 25-29, 1-18, 14-15, 12-22, 8-21, and 9-15 are 1000, 1200, 1300, 1300, 1400, 1400, 1400, and 1400 \$/period, respectively. The penalty cost for loss of load is 35000 \$/kW in each period.

\subsubsection{Estimated values of post-decision states}
\begin{figure}[!h]
	\centering
	\includegraphics[width=6cm]{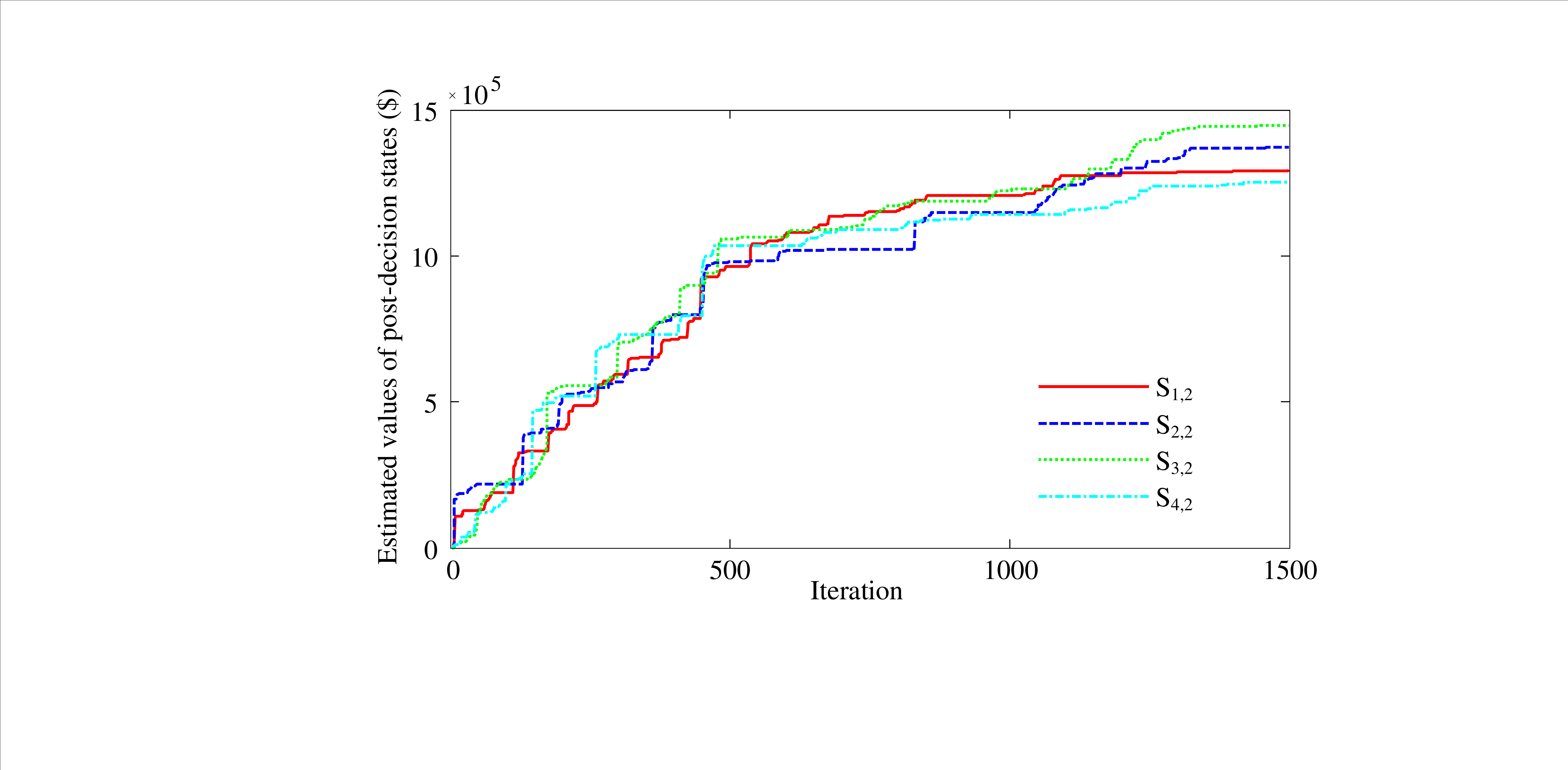}
	\caption{Iterations for estimated values of post-decision states}
	\label{figure1_fig}
\end{figure}

Because post-decision states are introduced to make the proposed recursively state-based model easily to be solved, one important task is to first estimate the values of these post-decision states according to the ADP algorithm. Due to a large number of post-decision states, we only show the estimated values of some post-decision states for the sake of exposition. Fig. \ref{figure1_fig} shows the estimated values of four post-decision states $S_{1,2}$, $S_{2,2}$, $S_{3,2}$ and $S_{4,2}$, shown in Table \ref{TableStateEx}, in the second decision period. 1500 iterations were performed to get the estimated values of post-decision states, and the estimated values of the post-decision states $S_{1,2}$, $S_{2,2}$, $S_{3,2}$ and $S_{4,2}$ converge to $1.29 \times 10^6 \$$, $1.37 \times 10^6 \$$, $1.45 \times 10^6 \$$, and $1.24 \times 10^6 \$$, respectively. With the estimated values, the term $v_t^{a_t}(S_{i,t}^{a_t})$ is known when optimizating (\ref{New_GenericModel1}). In this case, the stochastic problem is transformed into a deterministic problem. 
\begin{table}[H]
	\centering
	\renewcommand{\arraystretch}{1.3}
	\caption{Post-decision states $S_{1,2}$, $S_{2,2}$, $S_{3,2}$ and $S_{4,2}$ }
	\label{TableStateEx}
	\begin{tabular}{p{1.45cm}<{\centering}p{2.45cm}<{\centering}p{3.65cm}<{\centering}}
		\hline
		Time Period    & Post-decision states    & Disconnected lines   \\ \hline
		$2$    & $S_{1,2}$               & 10-11, 12-13, 25-29, 1-18, 8-21     \\ 
		$2$    & $S_{2,2}$               & 10-11, 25-29, 1-18, 14-15, 8-21     \\ 
		$2$    & $S_{3,2}$               & 10-11, 25-29, 1-18, 12-22, 8-21     \\ 
		$2$    & $S_{4,2}$               & 25-29, 1-18, 14-15, 12-22, 8-21     \\
		$1$    & $S_{1,1}$               & 10-11, 12-13, 25-29, 1-18, 8-21     \\ \hline
	\end{tabular}
\end{table} 

\begin{figure}[!h]
	\centering
	\includegraphics[width=6cm]{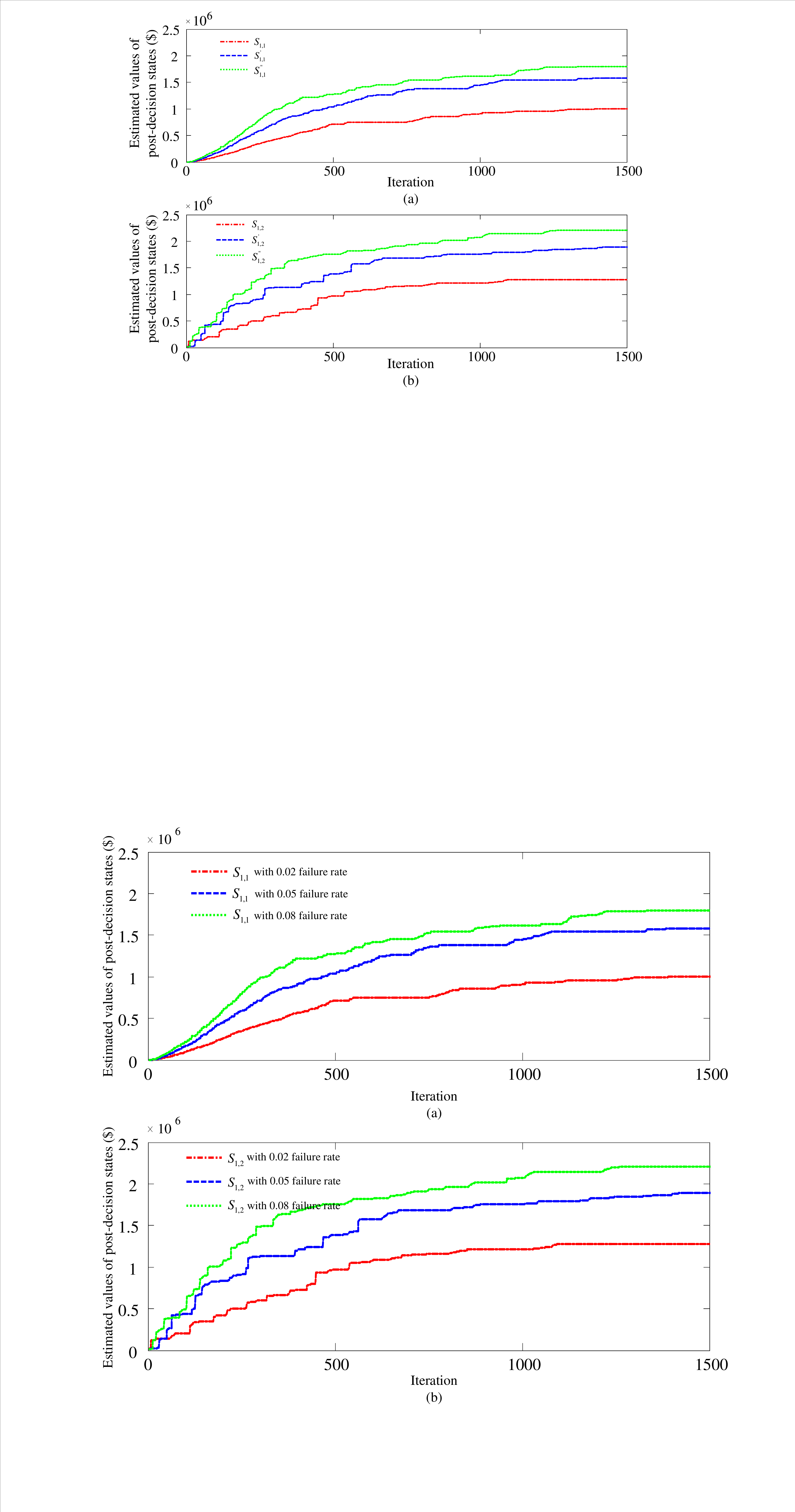}
	\caption{Iterations for estimated values of post-decision states with different failure rates}
	\label{figure2_fig}
\end{figure}

Different intensity of severe weather will result in different failure rates, which have great impacts on dispatch strategies. From the perspective of the mathematical model, different failure rates cause different estimated values of each post-decision state. Fig. \ref{figure2_fig} shows the estimated values of several post-decision states with different failure rates. For example, the estimated values of the post-decision state $S_{1,1}$ in the first period are $1.01 \times 10^6 \$$, $1.58 \times 10^6 \$$, and $1.81 \times 10^6 \$$ when the failure rates are 0.02, 0.04, and 0.06, respectively. %The estimated values of the post-decision state $S_{1,2}$ in the second period are $1.27 \times 10^6 \$$, $1.89 \times 10^6 \$$, and $2.19 \times 10^6 \$$ when the failure rates are $0.02$, $0.04$, and $0.06$, respectively. It is observed that a higher failure rate cause a larger estimated value of a post-decision state.  

\begin{figure}[!h]
	\centering
	\includegraphics[width=6cm]{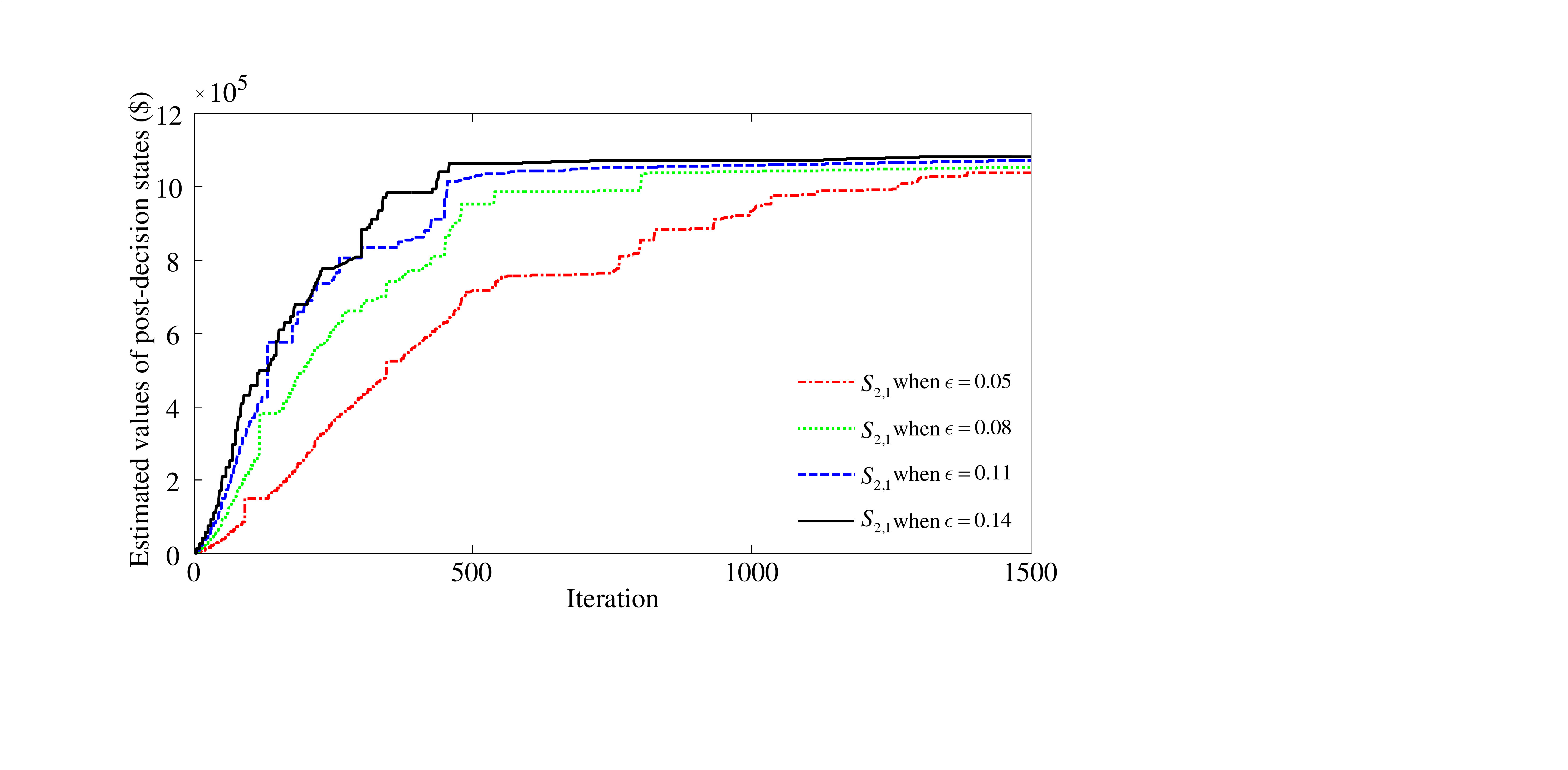}
	\caption{Iterations for estimated values of post-decision states with different values of $\epsilon$}
	\label{figure3_fig}
\end{figure}

When updating the estimated values of post-decision states by using (\ref{Update_model}), $\epsilon$ is artificially set. Fig. \ref{figure3_fig} shows the impacts of different values of $\epsilon$ on the estimated values of the post-decision state $S_{2,1}$. It is observed that the estimated value are close even when $\epsilon$ has different values.       


\subsubsection{Dispatch strategies with estimated values of post-decision states}
With the estimated values of each post-decision state based on ADP, the strategy corresponding to one observed real-time state can be obtained by a one-period deterministic optimization problem. Table \ref{IEEE33_Str1} and Table \ref{IEEE33_Str2} show the state-based strategies, and the original topology has the disconnected lines 10-11, 8-21, 9-15, 1-18, and 25-29. It is observed that the strategies make that the feeders impacted by the typhoon are downstream. This is reasonable because downstream feeders cause smaller outages even they are in failure due to the typhoon. %When the system operators do nothing during the unfolding event, the values of post-decision states $S_{1,2}$, $S_{2,2}$, $S_{3,2}$ and $S_{4,2}$ are $2.13 \times 10^6 \$$, $2.56 \times 10^6 \$$, $2.71 \times 10^6 \$$, and $2.01 \times 10^6 \$$, respectively. They are larger than the values $1.27 \times 10^6 \$$, $1.89 \times 10^6 \$$, and $2.19 \times 10^6 \$$ by using the proposed method, respectively. Because the values of post-decision states are expected costs in the subsequent periods, smaller values indicates better strategies.    

\begin{table}[H]
	\renewcommand{\arraystretch}{1.3}
	\caption{The first case of state-based Strategy}
	\label{IEEE33_Str1}
\begin{tabular}{cccc}
	\hline
	\multirow{2}{*}{Time Period} & \multirow{2}{*}{\begin{tabular}[c]{@{}c@{}}Observed State\\ (Component failure)\end{tabular}} & \multicolumn{2}{c}{Strategy}                                                                                                    \\ \cline{3-4} 
	&                                    & Open lines                    & Close lines                 \\ \hline
	1                 & -                & 19-20                                                            & 10-11                                                             \\ 
	2                 & 6-7              & 10-11, 6-26                                                      & 19-20, 8-21, 1-18                                                             \\ 
	3                 & -                & -                                                                & -                                                        \\ 
	4                 & -                & -                                                                & -                                                             \\ 
	5                 & -                & -                                                                & -                                                             \\ 
	6                 & -                & -                                                                & - \\  \hline
\end{tabular}
\end{table}

\begin{table}[H]
	\renewcommand{\arraystretch}{1.3}
	\caption{The second case of state-based Strategy}
	\label{IEEE33_Str2}
\begin{tabular}{cccc}
	\hline
	\multirow{2}{*}{Time Period} & \multirow{2}{*}{\begin{tabular}[c]{@{}c@{}}Observed State\\ (Component failure)\end{tabular}} & \multicolumn{2}{c}{Strategy}                                                                                                    \\ \cline{3-4} 
	&                                    & Open lines                    & Close lines                 \\ \hline
	1                 & -                & 19-20                                                            & 10-11                                                             \\ 
	2                 & -                & 6-26                                                      & 19-20, 1-18                                                             \\ 
	3                 & -                & -                                                                & -                                                        \\ 
	4                 & -                & -                                                                & -                                                             \\ 
	5                 & -                & -                                                                & -                                                             \\ 
	6                 & -                & -                                                                & - \\  \hline
\end{tabular}
\end{table}

\subsection{IEEE 123-bus system}
\subsubsection{Data description}
\begin{figure}[!h]
	\centering
	\includegraphics[width=6cm]{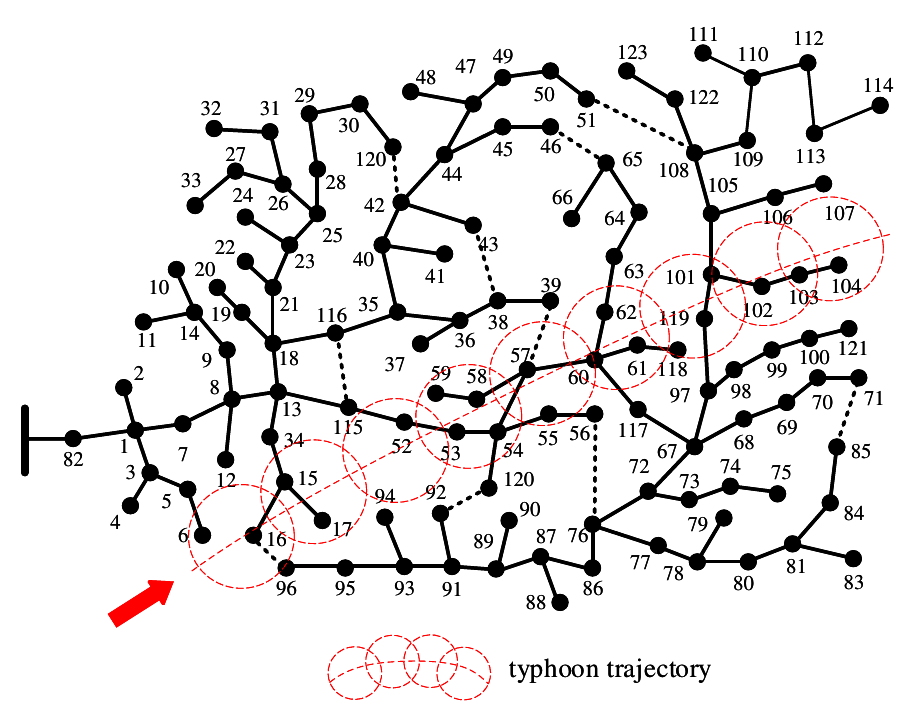}
	\caption{Topology of IEEE 123-bus system}
	\label{IEEE123_fig}
\end{figure}

Fig. \ref{IEEE123_fig} shows the topology of the IEEE 123-bus system and the trajectory  of a typhoon. The original topology has the disconnected lines 16-96, 92-120, 115-116, 42-120, 38-43, 39-57, 56-76, 46-65, 51-108, and 71-85. The lines 16-96, 92-120, 56-76, 39-57, 38-43, 42-120, 46-65, 51-108, 71-85, 52-53, 60-57, 60-117, 101-119, 63-64, and 67-117 are dispatchable. %The operational costs for the dispatchable lines are 1200 \$/period. The penalty cost for loss of load is 35000 \$/kW in each period.  

\subsubsection{Simulated results} 
\begin{figure}[!h]
	\centering
	\includegraphics[width=9.2cm]{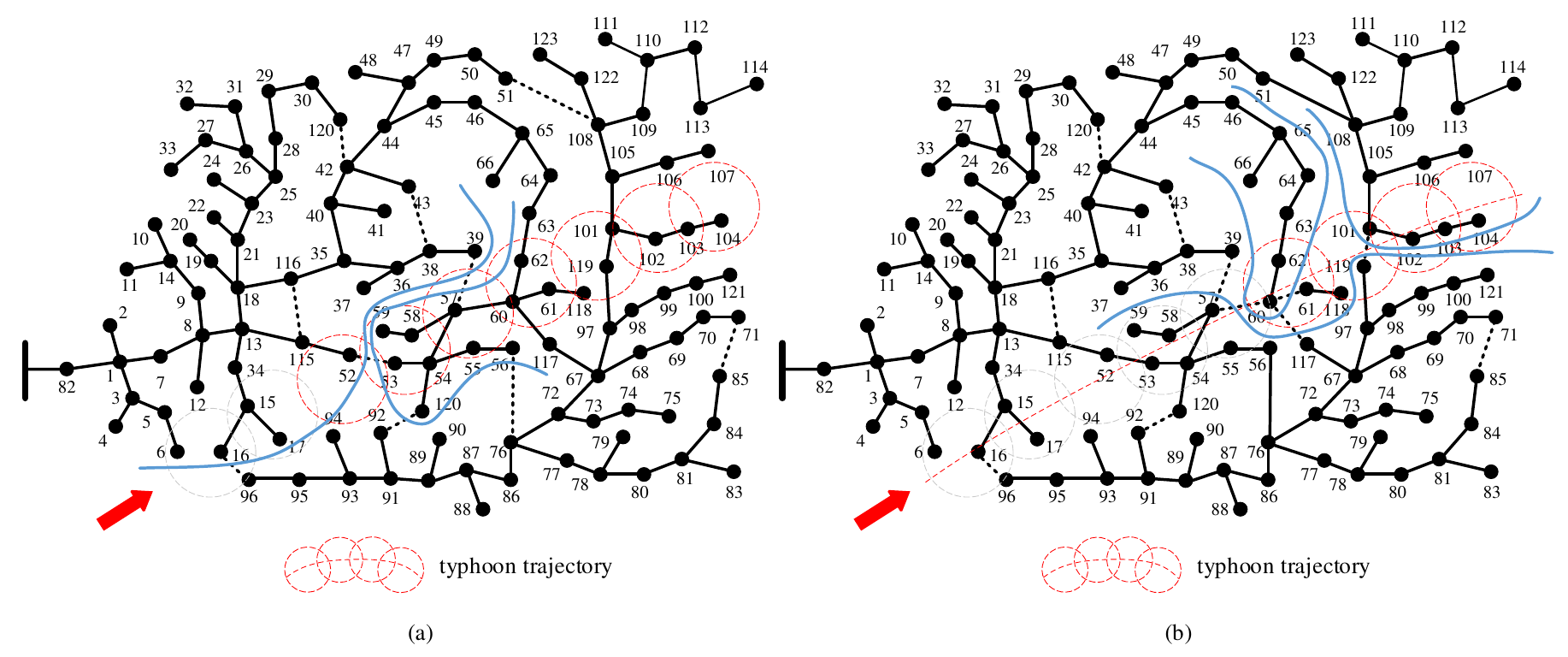}
	\caption{System Topologies (a) in the $3^{rd}$ period and (b) in the $6^{th}$ period.}
	\label{IEEE123_Strategy_fig}
\end{figure}

Based on the ADP algorithm, the estimated values of post-decision states can be obtained, and then the state-based strategies can be optimized. Table \ref{IEEE123_Str1} shows the state-based strategies on the trajectory of the typhoon, and the states on the trajectory are assumed to be generated stochastically based on failure rates caused by the typhoon. Fig. \ref{IEEE123_Strategy_fig} (a) and (b) show the topologies after implementing the state-based strategies in the $3^{th}$ and  $6^{th}$ periods, respectively. In the $3^{th}$ period, the line 52-53 is disconnected and the line 46-65 is connected to avoid balck out of downstream feeders if the typhoon fails the line 52-53. In the $6^{th}$ period, three lines (57-60, 60-117, 101-119) are disconnected and three lines (52-53, 56-76, 51-108) are connected to reduce possible black-out areas. It is observed that the state-based strategies try to make the feeders on the trajectory locate the terminal of the whole network to reduce potential loss of load. % Define the topologies in Fig. \ref{IEEE123_Strategy_fig} (a) and (b) as the post-decision states $S_1$ and $S_2$. When the system operators do nothing during the unfolding event, the values of $S_1$ and $S_2$ are $2.79 \times 10^6 \$$ and $2.03 \times 10^6 \$$, respectively. When using the proposed method, the values of $S_1$ and $S_2$ are $2.04 \times 10^6 \$$ and $1.79 \times 10^6 \$$, respectively. Smaller values indicates better strategies.  

\begin{table}[!h]
	\renewcommand{\arraystretch}{1.15}
	\centering
	\caption{State-based Strategy for IEEE 123-bus system}
	\label{IEEE123_Str1}
	\begin{tabular}{cccc}
		\hline
		\multirow{2}{*}{Time Period} & \multirow{2}{*}{\begin{tabular}[c]{@{}c@{}}Observed State\\ (Component failure)\end{tabular}} & \multicolumn{2}{c}{Strategy}                                                                                                    \\ \cline{3-4} 
		&                                    & Open lines                    & Close lines                 \\ \hline
		1                 & -                & -                                                                & -                                                             \\ 
		2                 & 15-17            & -                                                                & -                                                             \\ 
		3                 & -                & 52-53                                                            & 46-65                                                         \\ 
		4                 & -                & -                                                                & -                                                             \\ 
		5                 & 58-57            & -                                                                & -                                                             \\ 
		6                 & -                & \begin{tabular}[c]{@{}c@{}}57-60, 60-117,\\ 101-119\end{tabular} & \begin{tabular}[c]{@{}c@{}}52-53, 56-76\\ 51-108\end{tabular} \\ 
		7                 & -                & -                                                                & -                                                             \\
		8                 & 102-103          & -                                                                & -                                                             \\ 
		9                 & -                & -                                                                & -                                                             \\ \hline
	\end{tabular}
\end{table}

\section{Conclusion}
This paper proposed a Markov state-based decision-making model with dispatching system topology to improve the distribution system resilience throughout the unfolding events. The sequentially states of system topologies changed by the unfolding events and actions are modeled as Markov states, and the uncertainties between different Markov states are represented as transition probabilities that are determined by the component failure rates caused by the unfolding events. Based on Markov states, a recursive optimization model based on Markov decision processes, including the current cost and the expected cost in the future, is developed to make state-based actions at each decision time. To deal with `curse of dimensionality' caused by uncertainties, an approximate dynamic programming (ADP) approach with post-decision states and iteration is employed to solve the proposed model. With the estimated values of post-decision states, the stochastic problem with sequential multi-period stochastic optimization problem is transformed into a one-period deterministic problem. Case studies demonstrate that the state-based strategies try to make the feeders on the trajectory locate the terminal of the whole network to reduce potential loss of load, and in consequence to improve system resilience.

\bibliographystyle{IEEEtran}
\bibliography{IEEEabrv,RefDatabase}

%\begin{IEEEbiographynophoto}{Name (S'10)} received .
%\end{IEEEbiographynophoto}

\end{document}